# Max-Semi-Selfdecomposable Laws and Related Processes

S Satheesh

NEELOLPALAM, S. N. Park Road
Trichur – 680 004, **India.**

ssatheesh@sancharnet.in

and

E Sandhya

Department of Statistics, Prajyoti Niketan College
Pudukkad, Trichur – 680 301, **India.**

esandhya@hotmail.com

**Abstract:** Methods of construction of Max-semi-selfdecompsable laws are given. Implications of this method in random time changed extremal processes are discussed. Max-autoregressive model is introduced and characterized using the max-semi-selfdecompsable laws and exponential max-semi-stable laws. Some comments regarding the selfdecomposability of max-semi-stable laws are given.

## 1. Introduction.

This investigation is motivated by the discussion of selfdecomposable (SD) and semi-selfdecomposable (semi-SD) laws and its relation to autoregressive (AR(1)) series and subordination of Levy processes in Satheesh and Sandhya (2004*b*). The class of semi-SD laws was characterized therein using an additive AR(1) model. That there is a close relationship between SD laws and AR (1) series had been observed by Gaver and Lewis (1980). Parallel to SD and semi-SD laws, Pancheva (1994) and Becker-Kern (2001) have discussed their maximum versions *viz.* max-SD and max-semi-SD laws. Again, paralleling subordination of Levy processes Pancheva, *et al.* (2003) has discussed a random time changed or compound extremal processes (EP). However, in the literature we have not come across a maximum version of the AR(1) model and so here we introduce such a model and characterize it using max-semi-SD laws.

EPs are processes with increasing right continuous sample paths and independent max-increments. The univariate marginals of an EP determine its finite dimensional distributions. In general EPs are discussed with state space in $\mathbf{R}^d$, $d>1$ integer. This is because their max-increments are max-infinitely divisible (max-ID) and a discussion of max-ID laws makes sense only for distribution functions (d.f) on $\mathbf{R}^d$, $d>1$ integer as all d.fs on $\mathbf{R}$ are max-ID. However, we can discuss max-stable, max-semistable, max-SD and max-semi-SD laws on $\mathbf{R}$. Selfsimilar (SS) processes are those that are invariant in distribution under suitable scaling of time and space. Pancheva (1998, 2000) has developed SS and semi-SS EPs and showed that EPs whose univariate marginals are strictly max-stable (strictly max-semistable) constitute a class of SS (semi-SS) EPs having homogeneous max-increments. She also considered EPs whose univariate marginals are max-SD (max-semi-SD). Pancheva (1994, 1998, 2000) has described these laws and related processes in more generality, in terms of the invariance of their d.fs w.r.t one-parameter groups and cyclic groups of time-space changes. Here we will restrict our discussion to d.fs on $\mathbf{R}$ and the invariance w.r.t linear normalization. In the maximum scheme we have the following definitions.

**Definition.1.1** (Megyesi, 2002). A non-degenerate d.f $F$ is max-semi-stable($a,b$) if

$$F(x) = exp\{-x^{-\alpha} h(ln(x))\}, \; x>0, \; \alpha>0, \tag{1}$$

where $h(x)$ is a positive bounded periodic function with period $ln(b)$, $b>1$, and there exists an $a>1$ such that $ab^{-\alpha} = 1$. This class is denoted by $\boldsymbol{\Phi}_{\alpha,a,b}$, and,

$$F(x) = exp\{-|x|^{\alpha} h(ln(|x|))\}, \; x<0, \; \alpha>0, \tag{2}$$

where $h(x)$ is a positive bounded periodic function with period $|ln(b)|$, $b<1$, and there exists an $a>1$ such that $ab^\alpha = 1$. This class is denoted by $\Psi_{\alpha,a,b}$.

$\Phi_{\alpha,a,b}$ is the extended Frechet type and $\Psi_{\alpha,a,b}$ the extended Weibull. The extended Gumbel has the translational invariance property and is not considered here.

**Definition.1.2** (Becker-Kern, 2001). A non-degenerate d.f $F$ is max-semi-SD($c$) if for some $c>1$ and $v\in\mathbf{R}$ there is a non-degenerate d.f H such that

$$F(x) = F(c^v x + \beta)\ H(x)\ \forall\ x\in\mathbf{R}, \tag{3}$$

where $\beta=0$ if $v\neq0$ and $\beta = ln(c)$ if $v=0$. If (3) holds for every $c>1$, then $F$ is max-SD.

Pancheva, *et al.* (2003) has discussed random time changed or compound EPs and their theorem.3.1 together with property.3.2 reads: Let $\{Y(t),\ t\geq0\}$ be an EP having homogeneous max-increments with d.f $F_t(y) = exp\{-t\mu([\ell,y]^c)\}$, $y\geq\ell$, $\ell$ being the bottom of the rectangle $\{F>0\}$ and $\mu$ the exponential measure of $Y(1)$, that is, $\mu([\ell,y]^c) = -ln(F(y))$. Let $\{T(t),\ t\geq0\}$ be a non-negative process independent of $Y(t)$ having stationary, independent and additive increments with Laplace transform (LT) $\varphi^t$. If $\{X(t),\ t\geq0\}$ is the compound EP obtained by randomizing the time parameter of $Y(t)$ by $T(t)$, then $X(t) = Y(T(t))$ and its d.f is:

$$P\{X(t)<\mathrm{x}\} = \{\varphi(\mu([\ell,y]^c))\}^t. \tag{4}$$

Pancheva, *et al.* (2003) also showed that in the above setup $Y(T(t))$ is also an EP. For a description of the process $\{T(t)\}$ in this paper we need the following notion also.

**Definition.1.3** (Maejima an Naito, 1998). A probability distribution with characteristic function (CF) $f$ is semi-SD($c$) if for some $0<c<1$ there exists a CF $f_o$ such that

$$f(s) = f(cs) f_0(s), \ \forall \ s \in \mathbf{R}.$$

If this relation holds for every $0 < c < 1$ then $f$ is SD.

A stochastically continuous process $\{X(t), \ t \geq 0\}$ having stationary, independent and additive increments and $X(o) = 0$ is called a Levy process. A Levy process $X(t)$ such that $X(1)$ is SD (semi-SD) will be called a SD (semi-SD) process. Since SD and semi-SD laws are ID, above Levy processes are well defined.

**Note.1** For brevity the above EPs will be referred to as max-stable, max-semi-stable, max-SD, max-semi-SD EPs and the Levy process above as SD and semi-SD processes.

**Note.2** In this paper stable and semi-stable laws, selfsimilarity and semi-selfsimilarity are considered in the strict sense only.

With this background we give methods to construct max-semi-SD laws in section.2. Implications of these constructions in the context of compound EPs are then discussed. In section.3 we introduce a max-AR(1) series and characterize it using the max-semi-SD law. The max-AR(1) model is then modified to characterize the exponential max-semi-stable laws. Finally in section.4, certain selfdecomposability properties of max-semi-stable laws are given.

**2. Max-semi-SD laws and compound EPs.**

We begin with some remarks.

**Remark.2.1** If $h(x)$ in (1) is periodic w.r.t $ln(b_1)$ and $ln(b_2)$ such that $ln(b_1)/ ln(b_2)$ is irrational then $h(x) = \lambda$, a constant. When $h(x)$ is a constant $F$ is max-stable.

**Remark.2.2** The d.fs in $\boldsymbol{\Phi}_{\alpha,a,b}$ can be represented in the form $exp\{-\psi(x)\}$, where $\psi(x)$ satisfies $\psi(x) = a\psi(bx)$, for some $a > 1$, $b > 1$, and $\alpha > 0$ satisfying $ab^{-\alpha} = 1$. Similarly those in $\Psi_{\alpha,a,b}$ can be represented as $exp\{-\psi(x)\}$, where $\psi(x) = a\psi(bx)$, for some $a>1$, $b<1$, and $\alpha>0$ satisfying $ab^{\alpha} = 1$. On the other hand, the general solution to $\psi(x) = a\psi(bx)$, is $\psi(x) = x^{-\alpha} h(ln(x))$, $x>0$ or $|x|^{\alpha} h(ln(|x|))$, $x<0$, where $h(x)$ is as in (1) which can be proved along the lines in Lin (1994) or Pillai and Anil (1996).

**Remark.2.3** An EP $\{Y(t)\}$ is SS if for any $b>0$ there is an exponent $H>0$ such that

$$\{Y(bt)\} \overset{d}{=} \{b^{H}Y(t)\}. \qquad (5)$$

If (5) holds for some $b>0$ only then $\{Y(t)\}$ is semi-SS. Without loss of generality we may consider the range of $b$ as $0<b<1$ in the description of selfsimilarity because (5) is equivalent to $\{(b^{-1})^{H}Y(t)\} \overset{d}{=} \{Y(b^{-1}t)\}$ and thus the whole range $b>0$ is covered.

**Remark.2.4** In definition.1.2 of max-semi-SD($c$) laws, the case $v=0$ is not considered here. Thus the combination $c>1$ and $v\in\mathbf{R}-\{0\}$ implies that we are considering

$$F(x) = F(cx) H(x) \ \forall \ x\in\mathbf{R} \text{ and for some } c\in(0,1)\cup(1,\infty). \qquad (6)$$

Thus there is a one-to-one correspondence between the possible scale changes in max-semi-stable and max-semi-SD laws. $F$ is max-SD if (6) holds for every $c\in(0,1)\cup(1,\infty)$.

Satheesh (2002) and Satheesh and Sandhya (2004$a$) have discussed φ-max-stable and φ-max-semi-stable laws as follows. For a LT $\varphi$ the d.f $\varphi\{-ln(F(x))\}$ is φ-max-stable (φ-max-semi-stable) if the d.f $F(x)$ is max-stable (max-semi-stable). When $\varphi$ is exponential we call the corresponding d.f as exponential max-stable (max-semi-stable). Satheesh (2002) has discussed the connection among; mixtures of max-ID laws (theorem.2.1$b$), randomizing the time parameter of an EP (theorem.2.2$b$), construction of SD laws (property.2.1) and the possibility of constructing its max-analogue using mixtures of max-stable laws. This analogue is:

**Theorem.2.1** φ-max-stable laws corresponding to (1) and (2) are max-SD if $\varphi$ is SD.

Generalizing this Satheesh and Sandhya (2004$a$) proved the following;

**Theorem.2.2** φ-max-semi-stable($a,b$) laws corresponding to (1) and (2) are max-semi-SD($b$) if $\varphi$ is SD.

We now generalize this as follows.

**Theorem.2.3** φ-max-semi-stable($a,b$) laws corresponding to (1) (respectively (2)) are max-semi-SD($b$) if $\varphi$ is semi-SD($b^{-\alpha}$), $b>1$ (respectively semi-SD($b^{\alpha}$), $b<1$).

*Proof.* If $\varphi$ is semi-SD($b^{-\alpha}$), $b>1$, then there exists a LT $\varphi_o$ and the d.f of a φ-max-semi-stable($a,b$) law is;

$$\varphi\{\psi(x)\} = \varphi\{b^{-\alpha}\psi(x)\}\varphi_o\{\psi(x)\}, \ \forall \ x \in \mathbf{R},$$

Since $\psi(x) = a\psi(bx)$, for $a>1$, $b>1$, and $\alpha>0$ satisfying $ab^{-\alpha} = 1$, we have;

$$\varphi\{\psi(x)\} = \varphi\{\psi(bx)\}\varphi_o\{\psi(x)\} \text{ for some } b>1,$$

which completes the proof corresponding to the case (1). For the case corresponding to (2) we can proceed along the above lines and we have d.fs satisfying

$$\varphi\{\psi(x)\} = \varphi\{\psi(bx)\}\varphi_o\{\psi(x)\} \text{ for some } b<1.$$

This completes the proof of the theorem.

**Example.1** Consider the d.f $F(x) = \{1+\psi(x)\}^{-\beta}$, where $\psi$ is as in remark.2.2 and $\beta>0$. Since the gamma(1,$\beta$) law is SD the above d.f is max-semi-SD($b$) being that of a gamma-max-semi-stable($a,b$) law. By a similar line of argument, Satheesh and Sandhya (2004*b*) have shown that the CF $\{1+\psi(t)\}^{-\beta}$ is semi-SD(*b*), *b*<1. Write $b=\lambda^\alpha$ for $\lambda$<1 and $\alpha\in(0,2]$. Setting $\varphi$ to be this semi-SD($\lambda^\alpha$) the corresponding $\varphi$-max-semi-stable($a,\lambda$) laws are max-semi-SD($\lambda$).

In the case of compound EPs, extending the discussion in Pancheva, *et al.* (2003) we now have the following results. These results follow from the equation (4) and respectively from theorems 2.1, 2.2 and 2.3.

**Theorem.2.4** The EP obtained by compounding a max-stable EP is max-SD if the compounding process is SD.

**Theorem.2.5** The EP obtained from a random time changed max-semi-stable($a,b$) EP is max-semi-SD($b$) if the compounding process is SD.

**Theorem.2.6** The EP obtained by compounding a max-semi-stable($a,b$) EP corresponding to (1) (respectively (2)) is max-semi-SD($b$) if the compounding process is semi-SD($b^{-\alpha}$), $b$>1 (respectively semi-SD($b^\alpha$), $b$<1).

Now a condition for a semi-SS EP to be SS. The proof follows from remark.1.

**Theorem.2.7** If a max-semi-stable($a,b$) EP (hence semi-SS) $\{Y(t)\}$ satisfies (5) for two values $b_1$ and $b_2$ such that $ln(b_1)/ ln(b_2)$ is irrational then it is max-stable (hence SS).

## 3. Autoregressive models with a maximum structure.

We now develop AR(1) models with a maximum structure and give stationary solutions to it characterizing max-semi-SD and exponential max-semi-stable laws.

**Definition.3.1** A sequence $\{X_n\}$ of r.vs generates a max-AR(1) series if for some $\rho$>0 there exists an innovation sequence $\{\varepsilon_n\}$ of i.i.d r.vs such that

$X_n = \rho X_{n-1} \vee \varepsilon_n$ , $\forall$ $n > 0$ integer. $\qquad\qquad$ (7)

Since $X_{n-1}$ is a function only of $\varepsilon_j$, $j = 1, 2, \ldots, n-1$, it is independent of $\varepsilon_n$. Hence in terms of d.fs this is equivalent to;

$F_n(x) = F_{n-1}(x/\rho)\ F_\varepsilon(x)$ $\forall$ $x \in \mathbf{R}$.

Assuming the series to be marginally stationary we have;

$F(x) = F(x/\rho)\ F_\varepsilon(x)$ $\qquad\qquad$ (8)

Now comparing (8), (6) and (3) the following theorem is clear.

**Theorem.3.1** A sequence $\{X_n\}$ of r.vs generates a max-AR(1) series that is marginally stationary if and only if the distribution of $X_n$ is max-semi-SD($c$), $c = 1/\rho$.

**Remark.3.1** Unlike the additive AR(1) scheme where semi-SD($\rho$) laws characterizes only the case of $\rho < 1$, max-semi-SD($c$) laws characterizes the max-AR(1) scheme including the explosive situation of $\rho > 1$ as well, where $c = 1/\rho$.

Now let us modify the max-AR(1) scheme (7) as follows.

$X_n\ = \rho X_{n-1}$ , with probability $p$

$\quad = \rho X_{n-1} \vee \varepsilon_n$ , with probability $(1-p)$. $\qquad\qquad$ (9)

In terms of d.fs and assuming marginal stationarity of $X_n$ this is equivalent to;

$F(x) = F(x/\rho)\{p + (1-p)\ F_\varepsilon(x)\}$.

Further assuming $X_n$ and $\varepsilon_n$ to be identically distributed and setting $c = 1/\rho$ we have;

$F(x) = pF(cx)/\{1 - (1-p)F(cx)\}$. $\qquad\qquad$ (10)

This means that $F(x)$ is a geometric (geometric with mean $1/p$) maximum of its own type. Writing $F(x) = 1/\{1+\psi(x)\}$, we have from (10),

$$\frac{1}{1+\psi(x)} = \frac{p \frac{1}{1+\psi(cx)}}{1-(1-p)\frac{1}{1+\psi(cx)}} = \frac{1}{1+\frac{1}{p}\psi(cx)}, \forall\, x \in \mathbf{R.}$$

Hence, $\psi(x) = \frac{1}{p}\psi(cx)$, and thus $\psi(x)$ is now as in remark.2.2 with $a=\frac{1}{p}$ and $b=c=1/\rho$. Hence $F(x)$ is an exponential max-semi-stable($a$,$b$) law (since here $\varphi$ is the LT of an exponential law). Notice that here also both $\rho<1$ and $\rho>1$ are permitted. As the converse is clear we have proved;

**Theorem.3.2** A sequence $\{X_n\}$ of r.vs generates a max-AR(1) series with the structure of (9) that is marginally stationary and $X_n \overset{d}{=} \varepsilon_n$ if and only if the distribution of $X_n$ is an exponential max-semi-stable($a$,$b$) law, $a=\frac{1}{p}$ and $b=c=1/\rho$.

## 4. Selfdecomposability properties of max-semi-stable laws.

Satheesh and Sandhya (2004$b$) have shown that a semi-stable law is infinitely divisible and semi-SD. Here we show that similar results hold good for max-semi-stable laws as well.

**Theorem.4.1** Max-semi-stable($a$,$b$) laws corresponding to (1) and (2) are max-semi-SD($b$).

*Proof.* If $F(x)$ is a max-semi-stable($a$,$b$) d.f corresponding to (1) then we have;

$F(x) = \{F(bx)\}^a$, $x>0$, for some $a>1$, $b>1$, where $ab^{-\alpha} = 1$ for some $\alpha>0$.

$= F(bx)\{F(bx)\}^{a-1}$.

Hence by definition.1.2 $F$ is max-semi-SD($b$). A similar line of argument shows the fact for the case corresponding to (2). Notice that in the one-dimensional case $F$ is always max-ID.

**Theorem.4.2** Max-stable laws corresponding to (1) and (2) are max-SD.

*Proof.* If $F$ is a max-stable d.f then for each $b \in (0,1) \cup (1,\infty)$ there exists an $a > 1$ such that;

$$F(x) = \{F(bx)\}^a, x \in \mathbf{R}.$$

$$= F(bx)\{F(bx)\}^{a-1}, \text{ for every } b \in (0,1) \cup (1,\infty).$$

Hence $F$ is max-SD.

Thus we know how to construct max-semi-SD laws, its implication in compound EPs, related divisibility properties and the structure of max-AR(1) series.